\documentclass{article}
\usepackage{amssymb}
\usepackage{amsmath}
\usepackage{amsfonts}
\usepackage{graphicx}
\usepackage{epstopdf}
\usepackage{placeins}
\usepackage{color}
\usepackage{float}
\usepackage{caption}
\usepackage{enumitem}
\usepackage{tkz-euclide}
\usepackage{amsthm}
\usepackage[markup=underlined]{changes}

\providecommand{\keywords}[1]
{
  \small	
  \textbf{\textit{Keywords---}} #1
}

\title{{\bf {Uniqueness of universal dimensions and configurations of points and lines}\vspace{.2cm}}
\author{{\bf M.Y. Avetisyan\footnote{maneh.avetisyan@gmail.com, https://orcid.org/0000-0002-2060-2805} and R.L. Mkrtchyan\footnote{mrl55@list.ru, https://orcid.org/0000-0002-7565-3722 }}}
\date{}
}
\begin{document}
\maketitle
\date{ }
\begin{center}
  {\small {\it Yerevan Physics Institute, Yerevan, Armenia}}\\
\end{center}

\vspace{0.125cm}

\begin{abstract}
	
The problem of uniqueness of universal formulae for (quantum) dimensions of simple Lie algebras is investigated and 
connection of some of these functions with geometrical configurations, such as the famous Pappus-Brianchon-Pascal $(9_3)_1$ configuration 
of points and lines, is established by proposing some generic functions, which multiplied by a universal (quantum) dimension formula, preserve both its structure
and its values at the points from Vogel's table.
We deduce, that the appropriate realizable configuration $(144_3 36_{12})$ (yet to be found) will provide a symmetric non-uniqueness factor for any universal dimension formula.
\end{abstract}


\section*{Introduction}
 
Universal dimension formulae \cite{V0,V,Del,DM,Cohen,W3,LM1,MV,M16QD,AM,X2kng} are functions $f(\alpha,\beta,\gamma)$ of three  homogeneous so-called universal, or Vogel's parameters $(\alpha:\beta:\gamma)$ - the projective coordinates in Vogel's plane\footnote{In \cite{LM1} Vogel's plane is defined as the factor of 
projective plane by permutations of the homogeneous parameters. Here we refer
 as {\it Vogel's} the projective plane itself.}.
These universal parameters particularly parametrize all simple Lie algebras, by setting up
a correspondence between specific projective coordinates and the latters.
This correspondence is presented in Vogel's Table \ref{tab:V2}.
Universal formulae possess a specific structure: they are rational functions, where both the numerator and denominator decompose into products of the same number  of 
linear factors of universal parameters:

\begin{eqnarray}\label{Qx}
	F=\prod_{i=1}^{k}\frac{n_i\alpha+x_i\beta+y_i\gamma}{m_i\alpha+z_i\beta+t_i\gamma}
\end{eqnarray}

The general form of universal quantum dimensions is similar, with the linear functions substituted by their hyperbolic sines\footnote{ The following notation is used
	
	\begin{eqnarray*}
		\sinh\left[x: \right. \,\, \frac{A\cdot B...}{M\cdot N...}\equiv \frac{\sinh(xA)\sinh(xB)...}{\sinh(xM)\sinh(xN)...}
	\end{eqnarray*}
	}:

\begin{eqnarray}\label{Qxq}
	F(x)=\sinh[x: \prod_{i=1}^{k}\frac{n_i\alpha+x_i\beta+y_i\gamma}{m_i\alpha+z_i\beta+t_i\gamma} 
\end{eqnarray}

The characteristic property of a universal (quantum) dimension function is that at the points in the Vogel's
 projective plane, corresponding to each of the Lie algebras, it yields 
(quantum) dimension of some irreducible representation of the particular algebra, extended by the automorphism 
group of its Dynkin diagram \cite{Del,Cohen}.  
We also require that at the points, obtained by permutation of the universal coordinates for the initial ones, 
the output of the universal function is equal to a (quantum) dimension of some irreducible representation of the 
corresponding simple Lie algebra as well.
 It is permissible, that at some of the abovementioned points a universal function may be singular (see, however, \cite{X2kng,AMLR}, where the linear resolvability feature, which allows to assign relevant values to universal dimension formulae at their singular points, is described). 

As an example, we present the simplest universal quantum dimension, that is for the adjoint representation $ad$:

 \begin{eqnarray}  \label{cad}
	f(x) &=& -\frac{\text{sinh}\left(\frac{\gamma+2\beta+2\alpha}{4}x\right)}{\text{sinh}\left(\frac{\alpha}{4}x\right)}\frac{\text{sinh}
\left(\frac{2\gamma+\beta+2\alpha}{4}x\right)}{\text{sinh}\left(\frac{\beta}{4}x\right)}\frac{\text{sinh}\left(\frac{2\gamma+2\beta+\alpha}{4}x\right)}{\text{sinh}\left(\frac{\gamma}{4}x\right)}
\end{eqnarray}

We also present the universal quantum dimension $X(x,k,n,\alpha,\beta,\gamma)$
 of the Cartan products of the $k$-th power of $X_2$ and the $n$-th power of
  the adjoint representations \cite{X2kng}:

\begin{multline}
	X(x,k,n,\alpha,\beta,\gamma)=\\ 
	\sinh\left[\frac{x}{4}: \right.\prod _{i=0}^{k-1} \frac{(\alpha  (i-2)-2 \beta )^2(\alpha  (i-2)-2 \gamma )^2(\beta +\gamma +\alpha  (-(i-2))^2}{(\alpha  (i+1) )^2
		(\beta -\alpha  (i-1))^2(\gamma-\alpha  (i-1))^2}\times \\
	\times \prod _{i=0}^{n} \frac{(\alpha  (i+k-2)-2 \beta )(\alpha  (i+k-2)-2 \gamma )(\beta +\gamma +\alpha  (-(i+k-2))}
	{(\alpha (i+k+1))(\beta -\alpha  (i+k-1))(\gamma -\alpha(i+k-1))} \times \\
	\times \prod _{i=1}^{2k+n} \frac{(-\beta -2 \gamma +\alpha  (i-3))(-2 \beta -\gamma +\alpha  (i-3))(\alpha 
		(i-5)-2 (\beta +\gamma ))}
	{(\alpha  (i-2)-2 \beta )(\alpha  (i-2)-2 \gamma )(\beta +\gamma -\alpha  (i-2))} \times \\
	\times \frac{(\alpha +\beta)(\alpha +\gamma)(\alpha  (n+1))}{(2 \alpha +2 \beta )(2 \alpha +2 \gamma )(2 \alpha +\beta
		+\gamma )}\times \\
	\times \frac{(\alpha  (3 k+n-4)-2 (\beta +\gamma))(\alpha  (3 k+2 n-3)-2
		(\beta +\gamma ))}{(3 \alpha +2 \beta +2 \gamma )(4 \alpha +2 \beta +2
		\gamma )}
\end{multline}

Here $X_2$ is the representation, appearing in the following universal decomposition

\begin{eqnarray}
	\wedge^2 ad =ad \oplus X_2
\end{eqnarray}
which holds for all the simple Lie algebras.

As we see, universal formulae may be sufficiently intricate. One may ask
then, whether it is possible to rewrite them more compactly, without
violating their structure and preserving their values at the points from Vogel's table.

The investigation, presented in this paper, is a step forward to the answer to this question. 
We examine the problem of uniqueness of
 the universal dimension (quantum dimension) formulae.
 
Let $F_1$ and $F_2$ are universal (quantum) dimension formulae, yielding the same outputs at the points from Vogel's
 table, as well as at the corresponding points with permuted coordinates. 
Then their ratio $Q$ - {\it the non-uniqueness factor} is obviously equal to 1 at those points.

\begin{eqnarray}
	Q=\frac{F_1}{F_2}
\end{eqnarray}

Evidently, $Q$ has the same structure as (\ref{Qx}) or (\ref{Qxq}). 

So, the problem of uniqueness rephrases into the investigation of existence of 
such a $Q$, which is equal to $1$ both at the points from Vogel's table \ref{tab:V2}, 
as well as at all those, obtained by all possible permutations of the associated coordinates. 

The complete solution of this problem has not been achieved yet, however,

Below we propose expressions for several non-uniqueness factors which are equal to $1$ on particular distinguished
 lines in Vogel's plane, containing the points, related to the simple Lie algebras, which partially solves the problem of uniqueness. We then
  present a geometric interpretation of it, setting up a remarkable connection between the problem of uniqueness and the classical
study of configurations of points and lines.

First, note that the points from Vogel's table occupy the following distinguished lines:

\begin{eqnarray} \label{lines}
sl:	\alpha+\beta=0, \\ \label{sl}
so:	2\alpha+\beta=0, \\ \label{so}
sp:	\alpha+2\beta=0,\\
exc:	\gamma-2(\alpha+\beta)=0, \label{exc}
\end{eqnarray}
on which the linear, orthogonal, symplectic and the exceptional algebras are situated, 
respectively. 

We add an additional condition to the initial problem, namely,
 we require that $Q$ is equal to 1 not only at the points, associated
with the simple Lie algebras, but also on the entire distinguished lines. In case of the 
$sl,so,sp$ lines this requirement holds anyway, since a rational function, yielding $1$
 at the integer-valued coordinates on the halfline,
 is always equivalent to $1$ on the entire line. However, in case of the $exc$ line this 
 requirement is a non-trivial additional one, since it is a variation of Deligne's hypothesis \cite{Del,DM}, 
 which has not been proved yet. 

Next, the initial lines, as well as those, obtained by permutations (i.e. lines, obtained by permutations of the
 coordinates) collect into a set of $12$ lines, 
on which $Q$ is equivalent to $1$;
it includes the abovementioned $sl, so$ and $exc$ lines (we shall call them basic lines), and the following lines, 
obtained by all permutations, 
among which there is the $sp$ line as well (we omit $=0$ part in the lines' equations below):

\begin{eqnarray} \label{plines}\nonumber
		\alpha+\gamma, \beta+\gamma\\  \label{lines2}
	sp:	\alpha+2\beta, \alpha+2\gamma,\beta+2\gamma, \gamma+2\alpha, \gamma+2\beta\\ \nonumber
	\beta-2(\alpha+\gamma),\alpha-2(\beta+\gamma),
\end{eqnarray}

Our investigation is built in the following way: first, we require that $Q\equiv1$
 on each of the $3$ basic lines - the $sl, so, exc$, and derive a 
general non-uniqueness factor for this case in Section 1.1.

We then impose the requirement for $Q\equiv1$ on the $sp$ line -
one of the lines, obtained by the $so$ line after the $\alpha \leftrightarrow \beta$ permutation of the parameters,
and find several solutions, employing
results taken from the theory of configurations of points and lines in Section 2.

 The problem of further investigation of the $Q\equiv1$ requirement for the 
 remaining lines is formulated in Sections 2.5 and 2.6, both in the framework of the
 geometric and permutation invariance approach. 

\section{$Q$ for universal (quantum) dimensions}

To simplify calculations, we make the following change of coordinates :

\begin{tabular}{c|c}
\centering

$\alpha^{\prime}=\alpha+\beta$	& $\alpha=-\alpha^{\prime}+\beta^{\prime}$ \\
	\hline
$\beta^{\prime}= 2\alpha+\beta$	& $\beta=2\alpha^{\prime}-\beta^{\prime}$ \\
	\hline
	$\gamma^{\prime} = \gamma-2(\alpha+\beta) $	& $\gamma=2\alpha^{\prime}+\gamma^{\prime}$ \\
	
\end{tabular}

such that in the primed coordinates the equations of the basic lines $sl,so,exc$ will simply be 
$$\alpha^{\prime}=0, \beta^{\prime}=0,\gamma^{\prime}=0.$$ 
The line equations (\ref{plines}) in terms of these new coordinates 
rewrite as:

\begin{eqnarray}  \nonumber
	\alpha^{\prime}+\beta^{\prime}+\gamma^{\prime},  4\alpha^{\prime}-\beta^{\prime}+\gamma^{\prime}   \\
	  3\alpha^{\prime}-\beta^{\prime}, 3\alpha^{\prime}+\beta^{\prime}+2\gamma^{\prime}, 6\alpha^{\prime}-\beta^{\prime}+2\gamma^{\prime}, 2\beta^{\prime}+\gamma^{\prime}, 6\alpha^{\prime}-2\beta^{\prime}+\gamma^{\prime} \\ \nonumber	
	  -9\alpha^{\prime}+3\beta^{\prime}-2\gamma^{\prime},-3\beta^{\prime}-2\gamma^{\prime}
\end{eqnarray}

From now on we work with the primed parameters, dropping the prime mark for 
convenience (up to special notice).

\subsection{$Q$ for universal dimensions: $sl, so, exc$ and $sp$ lines}

Let's take the general form (\ref{Qx}) of $Q$, assuming it is written in terms of the primed parameters, and 
consider its values on the  three lines $\alpha=0, \beta=0, \gamma=0.$

We require $Q\equiv1$ at $\alpha=0$. 
Then 

\begin{eqnarray}\label{Qx2}
	1\equiv\prod_{i=1}^{k}\frac{x_i\beta+y_i\gamma}{z_i\beta+t_i\gamma}
\end{eqnarray}

and one deduces, that $z_i=l_ix_{q(i)}, t_i=l_iy_{q(i)}$, with some permutation $q(i), i=1,...k$, and non-zero multipliers $l_i$ with $l_1l_2...l_k=1$. 

Substituting these relations into (\ref{Qx}), one has 

\begin{eqnarray}\label{Qx22}
	Q=\prod_{i=1}^{k}\frac{n_i\alpha+x_i\beta+y_i\gamma}{m_i\alpha+l_ix_{q(i)}\beta+l_iy_{q(i)}\gamma}
\end{eqnarray}

Absorbing the $1/l_i$ into $m_i$, renumbering $m_i\rightarrow m_{q(i)}$ and changing the order of the multipliers in the denominator, we rewrite $Q$ as:

\begin{eqnarray}\label{Qx23}
	Q=\prod_{i=1}^{k}\frac{n_i\alpha+x_i\beta+y_i\gamma}{m_i\alpha+x_{i}\beta+y_{i}\gamma}
\end{eqnarray}

Now let $Q\equiv1$ at $\beta=0$:
\begin{eqnarray}\label{Qx24}
	1\equiv\prod_{i=1}^{k}\frac{n_i\alpha+y_i\gamma}{m_i\alpha+y_{i}\gamma}
\end{eqnarray}

Then one must have $y_i=k_iy_{s(i)}, m_i=k_in_{s(i)}$, with some permutation $s(i)$ and with $k_1k_2...k_k=1$, so that $Q$ accepts the form: 

\begin{eqnarray}\label{Qx25}
	Q=\prod_{i=1}^{k}\frac{n_i\alpha+x_i\beta+y_i\gamma}{k_i n_{s(i)}\alpha+x_{i}\beta+k_iy_{s(i)}\gamma}=
	\prod_{i=1}^{k}\frac{n_i\alpha+x_i\beta+y_i\gamma}{k_i n_{s(i)}\alpha+x_{i}\beta+y_i\gamma}
\end{eqnarray}

Next we require $Q\equiv1$ at $\gamma=0$:

\begin{eqnarray}\label{Qx26}
	1\equiv \prod_{i=1}^{k}\frac{n_i\alpha+x_i\beta}{ k_i n_{s(i)}\alpha+ x_i \beta}
\end{eqnarray}

Again, from this relation we infer  

\begin{eqnarray}
	x_i=c_ix_{p(i)} \\
	k_i n_{s(i)}=c_i n_{p(i)}
\end{eqnarray}

for some permutation $p(i)$ and $c_i$ with $c_1c_2...c_k=1$.

So, altogether we have the following expression for $Q$ with the restrictions on its parameters:

\begin{eqnarray}\label{Qx27}
	Q=\prod_{i=1}^{k}\frac{n_i\alpha+x_i\beta+y_i\gamma}{ k_i n_{s(i)}\alpha+ x_i\beta+y_{i}\gamma}=\prod_{i=1}^{k}\frac{ n_i\alpha+x_i\beta+y_i\gamma}{ c_i n_{p(i)}\alpha+ x_i\beta+y_{i}\gamma} \\    \label{3lineBeq}
	x_i=c_ix_{p(i)} \\ 
	y_i=k_iy_{s(i)} \\
	k_i n_{s(i)}=c_i n_{p(i)} \\
	c_1c_2...c_k=1 \\  \label{3lineEeq}
	k_1k_2...k_k=1 
\end{eqnarray}
for some permutations $s(i), p(i)$. Note that after having solved these equations, one must check the $Q$ on absence of any cancellation in it. 

It is easy to show, that there is not a non-trivial solution if $k=1,2$. 
For $k=3$ one can show that the existence of a non-trivial solution requires that the permutations $s(i),p(i)$ do not have fixed points and do not coincide,
 i.e. $s(i)=i+1, p(i)=i+2 \,\, (mod \,\,3)$, or vice versa. One can also show that $n_i \neq 0$, so that one can factor them out, 
 or effectively put $n_i=1$, so that $k_i=c_i$, $y_3=c_3 y_1, y_2=c_2 c_3 y_1, x_2=c_2 x_1, x_3=c_2 c_3 x_1$. Denoting $x_1=x, y_1=y$, we get the final expression of $Q$:

\begin{eqnarray} \label{Q33}
	\frac{(\alpha+\beta x+\gamma y) (\alpha c_1 c_2 +\beta c_2 x+\gamma y) (\alpha c_1+\beta c_1 c_2 x+\gamma y)}{(\alpha c_1+\beta x+\gamma y) (\alpha+\beta c_2 x+\gamma y) (\alpha c_1 c_2 +\beta c_1 c_2 x+\gamma y)}
\end{eqnarray}

Finally, we require $Q\equiv1$ when $3\alpha-\beta=0$:

\begin{eqnarray}
	1\equiv \prod_{i=1}^{k}\frac{\alpha(n_i+3x_i)+y_i\gamma}{\alpha(c_i n_{p(i)}+3 c_i x_i)+y_i\gamma}
\end{eqnarray}

which leads to
\begin{eqnarray}
	c_i n_{p(i)}+3x_i= r_i(n_{v(i)}+3x_{v(i)}) \\
	y_i=r_i y_{v(i)} \\
	\prod_{i=1}^k r_i=1 \\
	i=1,2,...,k
\end{eqnarray}
for some permutation $v(i)$. 

So, altogether we have the

{\bf Proposition 1.} {\it The general expression for a non-uniqueness factor $Q$ for universal dimensions, which is equal to $1$ on each of the
 $sl:\alpha=0, so:\beta=0$ and $exc:\gamma=0$ lines in Vogel's
plane, writes as follows:

\begin{eqnarray}\label{Qx27b}
	Q=\prod_{i=1}^{k}\frac{n_i\alpha+x_i\beta+y_i\gamma}{ k_i n_{s(i)}\alpha+ x_i\beta+y_{i}\gamma}=\prod_{i=1}^{k}\frac{ n_i\alpha+x_i\beta+y_i\gamma}{ c_i n_{p(i)}\alpha+ x_i\beta+y_{i}\gamma}  
\end{eqnarray}
	
with parameters, satisfying the following equations

\begin{eqnarray}	\label{4lineBeq}
 	x_i=c_ix_{p(i)} \\
	y_i=k_iy_{s(i)} \\   \label{end1} 
	k_i n_{s(i)}=c_i n_{p(i)} \\  \label{4-th-eq}
	y_i=r_i y_{v(i)} \\           \label{end2}
	c_i n_{p(i)}+3x_i= r_i(n_{v(i)}+3x_{v(i)}) \\
	c_1c_2...c_k=1 \\
	k_1k_2...k_k=1 \\   \label{4lineEeq}
	r_1 r_2...r_k=1 
\end{eqnarray}

for some permutations $s(i), p(i), v(i), \, i=1,2...k$.}

{\bf Remark.} As follows from the example above, one can get a trivial ($Q=1$) or a non-trivial non-uniqueness factor $Q$ depending on the particular choice of permutations. 

Evidently, one can keep on going in this manner line after line, requiring $Q\equiv1$ on each of the other $8$ lines. 
Note, that at each of the steps one has to choose a permutation - a cancellation pattern, 
as well as introduce the coefficients, such as $k_i$, satisfying
 $k_1k_2...k_k=1$.
 Then, provided these introduced coefficients are given,
  one has to solve the linear 
equations on the parameters. 
Finally, one should check whether the solution is non-trivial,
 i.e. make sure, that there is no cancellation in $Q$. 
 Since the number of permutations is $k!$, a lot of candidates for the solution will be obtained at each of the steps (each line).

We shall not try to derive and solve these equations directly, instead, we will derive a solution for $k=4$, 
after setting up a geometrical interpretation of 
the problem and using several results from the theory of configurations of points 
and lines in Section \ref{configur}.

\subsection{$Q$ for universal quantum dimensions}

{\bf Lemma}: 
{\it  If 
\begin{eqnarray}
\sinh[x:	\prod_{i=1}^{k}\frac{n_i}{m_i} \equiv1
\end{eqnarray}

then the $\{n_i, i=1,...,k\}$ and $\{\epsilon_i m_i, i=1,...,k\}$ sets coincide for some choice of $\epsilon_i=\pm 1, i=1,...k$, $\prod_{i=1}^{k}\epsilon_i=1$.}

\begin{proof}
Without loss of generality, let $\vert n_k\vert=max(\{\vert n_1\vert,...,\vert n_k\vert, \vert m_1 \vert,...,\vert m_k\vert \})$.
Then some of the factors in the numerator of $Q$ is zero at $x=i \pi / n_k$.
To cancel this zero out, there is necessarily a multiplier in the denominator, zeroing at the same point. Thus $i\pi /n_k=i\pi q/m_i$, for some $q\in \mathbb{Z}$, from which it follows that
$q=\pm 1$ and $m_i=\pm n_k$. Applying the same logic for the remaining factors, we prove the statement of the lemma.
\end{proof}

Now consider the general expression (\ref{Qxq}) of $Q$:

\begin{eqnarray}
	Q(x)=\sinh[x: \prod_{i=1}^{k}\frac{n_i\alpha+x_i\beta+y_i\gamma}{m_i\alpha+z_i\beta+t_i\gamma} 
\end{eqnarray}

Due to the Lemma the analysis of this expression is almost the same as that in the previous subsection. The only difference is that the parameters $k_i, c_i, r_i$, 
must now be equal $\pm 1$, so that one can write the following 

{\bf Proposition 2.}
{\it The general expression for a non-uniqueness factor $Q$ for universal quantum dimensions, which is equal to $1$ on each of the
 $sl:\alpha=0, so:\beta=0$ and $exc:\gamma=0$ lines in Vogel's
plane, writes as follows:

\begin{eqnarray}\label{Qx28}
	Q(x)=\sinh[x: \prod_{i=1}^{k}\frac{n_i\alpha+x_i\beta+y_i\gamma}{ k_i n_{s(i)}\alpha+ x_i\beta+y_{i}\gamma}=\sinh[x: \prod_{i=1}^{k}\frac{ n_i\alpha+x_i\beta+y_i\gamma}
	{ c_i n_{p(i)}\alpha+ x_i\beta+y_{i}\gamma} 
	\end{eqnarray}

with the same sets of equations: (\ref{3lineBeq}) - (\ref{3lineEeq}) for the case of three lines - $sl, so, exc$, and  (\ref{4lineBeq}) - (\ref{4lineEeq}) for the case of four lines 
$sl, so, sp, exc$, and with $k_i=\pm 1, c_i=\pm 1, r_i=\pm 1, \,\, i=1,...k$. }

The same Remark as in 1.1 applies to this case too. Actually it is easy to prove, that in the case of three lines, i.e. equations  (\ref{3lineBeq}) - 
(\ref{3lineEeq}) there is not a non-trivial solution $Q(x)$ for the quantum 
dimensions for $k$'s up to $k=3$, inclusively, 
so that the first non-trivial solution exists for $k=4$. 

A solution, obtained by employing some knowledge from the
theory of configurations, is presented below.

\section{Configurations of points and lines and the problem of uniqueness of universal formulae}
 \label{configur}

In this section we provide a geometric point of view to the problem
of uniqueness, setting up its connection with a classical problem of 
the so-called {\it configurations}, namely {\it configurations of points and lines}. 

\subsection{Geometric representation of universal formulae}
First, observe, that each of the linear factors in the expression of $Q$ corresponds to a line in the projective Vogel's plane. Indeed, to each of the factors $x\alpha+y\beta+z\gamma$ 
one can put in correspondence the line equation $x\alpha+y\beta+z\gamma=0$.

Thus, for any given expression for universal (quantum) dimension, with say $k$ multipliers, we can draw a unique picture in the Vogel's plane, consisting of $k$ lines, 
corresponding to the 
linear factors in numerator, which will be referred as {\it red lines} for convenience, and $k$ {\it green} lines for those in the denominator.
In addition, we can draw a number of {\it black} lines, 
corresponding to the distinguished $sl, so, exc$, lines as well as those, associated to the permuted coordinates - such as the $sp$ line.  

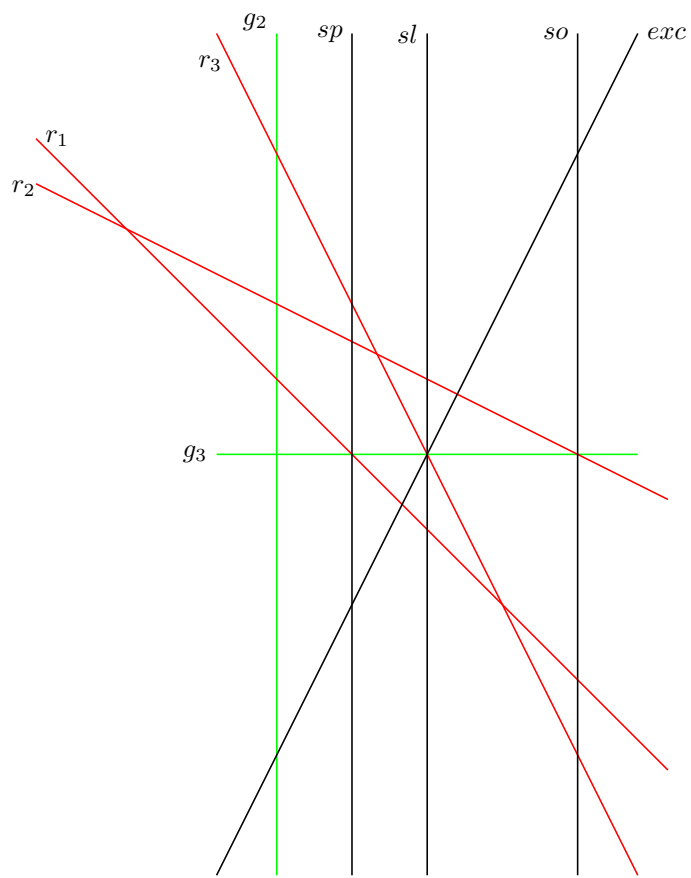
\begin{figure}
 \centering
\begin{tikzpicture} \tkzDefPoints{0/0/F, -2/3/A, 1/0/B, 2/0/C, 4/4/D, 4/0/E, 0/4/G, 0/-4/K, 4/-3/L, 4/-4/M, 1/4/N, 1/-4/J, 2/4/H, 2/-4/P}
\tkzDrawLines[color=green](K,G E,F)
\tkzDrawLines[color=red](E,A L,A G,M)
\tkzDrawLines(N,J H,P D,M K,D)
	\tkzLabelLine[pos=-0.22, left](G,K){$g_2$}
		\tkzLabelLine[pos=-0.2, left](F,E){$g_3$}
		\tkzLabelLine[pos=-0.18, left](A,E){$r_2$}
		\tkzLabelLine[pos=-0.2, right](A,L){$r_1$}
		\tkzLabelLine[pos=-0.15, left](G,M){$r_3$}
		\tkzLabelLine[pos=-0.2, left](N,J){$sp$}
		\tkzLabelLine[pos=-0.2, left](H,P){$sl$}
		\tkzLabelLine[pos=-0.2, left](D,M){$so$}
		\tkzLabelLine[pos=-0.2, right](D,K){$exc$}
\end{tikzpicture}
\caption{The "sketch" of the dimension formula (\ref{cad}).} \label{fig:Conf1}
\end{figure}

One can see the corresponding picture\footnote{The labeling of the lines in the following figures is meant to identify the corresponding colors
 they are given.
For example, in Figure \ref{fig:Conf1}, $r_2$ identifies the line, associated to the second factor in the numerator of (\ref{cad}), and $g_3$ - 
to the third factor in the corresponding denominator.} for the simplest universal formula, namely the dimension of the 
adjoint representation (\ref{cad}), in Figure \ref{fig:Conf1}.

Let's consider the picture, associated to a non-uniqueness factor $Q$.
 It turns out that each of the black lines must contain $k$ points, at which a green and a red line intersect. 

Indeed, this statement exactly rephrases the cancellation mechanism, described in the previous section:
 when restricting $Q$ to a black line, each of the factors from the numerator is proportional to some factor from the denominator.
This means that these two factors are zeroing simultaneously, meaning, that residing on a black and, say, a red line at once, we necessarily reside on a green line too. 

It is easy to notice, that this corresponding picture also contains information about the choice of the permutations $s(i), p(i), ...$, (see equations (\ref{Qx27b})-(\ref{4lineEeq})), -
the intersection points of three different-colored lines obviously define the pairs of cancelling factors, when restricting the function to each of the distinguished black lines. 

Thus, the picture of $k$ black, $k$ red and $k$ green lines, corresponding to a non-uniqueness factor $Q$, has the following characteristic feature:
 on each of the black lines there are $k$ points at which a red and a green line intersect. 
 Note that besides these points of intersection of three differently colored 
 lines, there may be some other intersection points, which however will not be 
 of interest for us. 

\subsection{Configurations}

Let's introduce the following standard definitions \cite{GKF,GB}:

{\bf Definition 1.}

{\it We say a line is incident with a point, (equivalently, a point is incident with a line) if it passes through it (equivalently, if it lies on it).}

{\bf Definition 2.}

{\it A { \bf configuration} $(p_{\gamma} l_{\pi})$ is a set of $p$ points and $l$ lines, such that every point is incident with precisely $\gamma$ of these lines and
every line is incident with precisely $\pi$ of these points.}

{\bf Remark 1.} Notice, that the total number of incidences, on one hand,
 is equal to $p \gamma$, and is $l \pi$, on the other hand, so that from 
Definition 2 it 
follows, that $p \gamma= l \pi$.

{\bf Remark 2.} If $p=l, \gamma=\pi$, the configuration is denoted by $(p_\gamma)$.

We see that the picture of $k$ black, $k$ red and $k$ green lines, possessing the feature described in the previous subsection,
 turns into a configuration iff the
number of black, red and green lines coincide and is equal to
 $k$.
Obviously, the corresponding configuration will be $(k^2_3,3k_k)$. 

\begin{figure}
 \centering
	\begin{tikzpicture} \tkzDefPoints{0/0/A, 6/0/B, 3/5/C, 
			1.4/2.33/D, 4.4/2.67/E, 3.3/0/F, 3.68/0.91/G, 3.41/2.56/H,
			2.05/1.54/I, 4.14/3.1/J, 1.13/1.89/K, 3.83/0/L}
\tkzDrawLines(A,B C,B A,C F,D D,E E,F A,J C,L B,K)
\tkzDrawPoints(A,B,C,D,E,F,G,H,I) 
			\end{tikzpicture}
	\caption{The $(9_3)_2$ configuration, which cannot be "colored" in order to be corresponded to some $Q$} \label{fig:Conf3}
\end{figure}
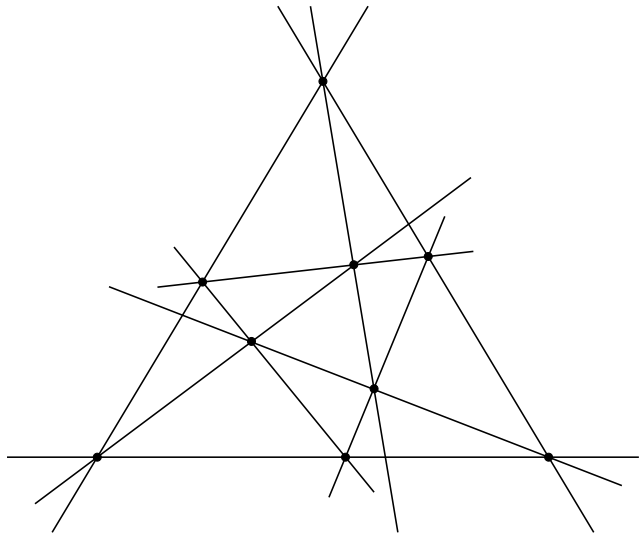

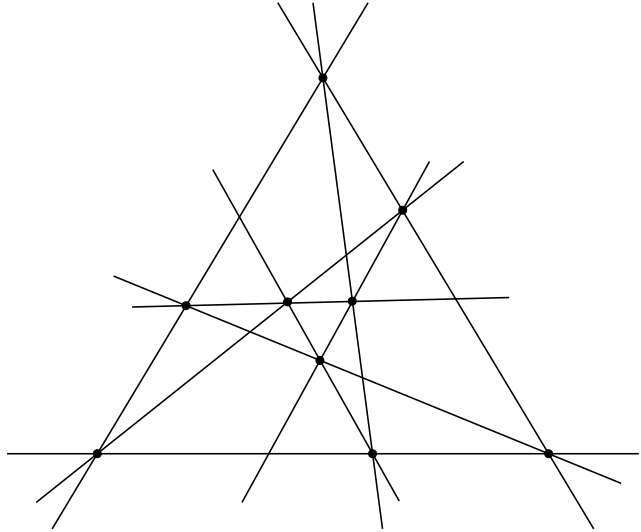
\begin{figure}
 \centering
\begin{tikzpicture} \tkzDefPoints{0/0/A, 6/0/B, 3/5/C, 3.66/0/D, 1.18/1.97/E, 4.06/3.24/F, 3.39/2.03/G, 2.96/1.24/H, 2.53/2.02/I, 1.89/3.15/J, 4.76/2.06/K, 2.28/0/L}
\tkzDrawLines(A,B C,B A,C A,F C,D B,E J,D E,K F,L)
 \tkzDrawPoints(A,B,C,D,E,F,H,I,G)
	\end{tikzpicture}
\caption{The $(9_3)_3$ uncolorable configuration} \label{fig:Conf4}
\end{figure}

However, if we have a configuration $(k^2_3,3k_k)$, it doesn't mean that we can definitely construct a corresponding $Q$. 
The possible obstacle is that one would not be able to attribute the black, red and green colors to its $3k$ lines such that at each of the points, 
belonging to the configuration, three lines of different colors meet. Such configuration are presented in Figures \ref{fig:Conf3} and \ref{fig:Conf4}.

For any given configuration $(p_{\gamma} l_{\pi})$ one can construct a so-called {\it configuration table}: we label the points and lines of that configuration, 
then for each of the lines allocate a column, consisting of the labels of the 
points, which are incident with the corresponding line. 
Characteristic properties of a configuration table are the following: the label for each of the points occurs in exactly $\gamma$ columns, different columns do not
contain two similar labels of points, and each column contains exactly $\pi$ labels. 
Two configuration tables are identical, if they coincide after some relabeling of points and lines, and/or rearranging the points in a given column.  

So, "possible" configurations of a given type $(p_{\gamma} l_{\pi})$ can be considered simply as different configuration tables of that type. 

Further, a configuration table is called {\it realizable}, if one can construct a geometrical picture of lines and points corresponding to it.
Not all tables are realizable.

\subsection{The $(9_3)_1$ configuration and $Q$ for the $sl, so, exc$ lines}

A relevant configuration happens to be corresponding to the solution (\ref{Q33}), derived for $k=3$, which is equivalent to $1$ on
 three basic lines  - (\ref{sl}), (\ref{so}), and (\ref{exc}).
It is the configuration $(9_3,9_3)_1$, which is usually referred as $(9_3)_1$ \cite{GKF,GB}, since the terms in the standard notation 
coincide. This configuration is also known as the Pappus (Pappus of Alexandria) or Pappus-Brianchon-Pascal configuration, 
which is presented in Figure \ref{fig:Conf2}.

The index in the notation $(9_3)_1$ is to indicate the fact, that there are several  $(9_3)$ configurations, so that it is used  to distinguish these.
Possible values of the index, i.e. the number of different configurations $(9_3)$ is $3$, equivalently, there are three different configuration tables for $(9_3)$ configuration. 
Each of these $3$ tables happens to be realizable. However, only one of them, presented in Figure \ref{fig:Conf2}, $(9_3)_1$ from \cite{GKF}, can be colored in the way we need. 
For example, for the configuration $(9_3)_2$ (see Figure \ref{fig:Conf3}) 
it is impossible to distinguish $3$ black lines, since for any two lines of the configuration, there is always a third one, which intersects with one of them at some point, which
 belongs to the configuration. This violates the requirement, that at each point of the configuration three lines of different
 colors intersect. The same reasoning holds for the remaining third configuration $(9_3)$ -  the $(9_3)_3$, see Figure \ref{fig:Conf4}, so that we 
 naturally can bring the following
 
 {\bf Proposition 3.}
{\it The non-uniqueness factor (\ref{Q33}) at $k=3$ is in one-to-one correspondence with
 the Pappus-Brianchon-Pascal $(9_3)_1$ configuration.}
 
 \begin{proof}
 The proof is obvious, since there is no other
  $(9_3)$ configuration which can be colored in the needed way. 
\end{proof}

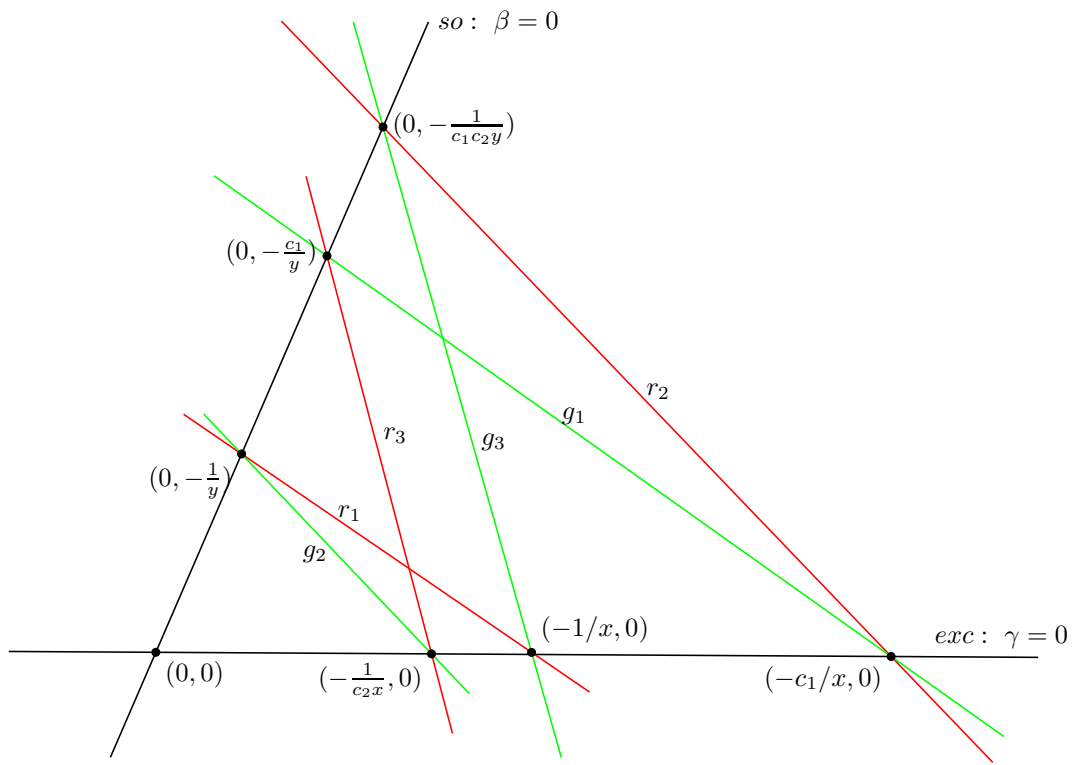
\begin{figure}
 \centering
\begin{tikzpicture} \tkzDefPoints{3.67/2.98/A, 5/3/B, 9.78/2.945/C, 
 1.145/5.64/H, 2.28/8.275/I, 3.025/9.99/J}
 \tkzInterLL(C,A)(J,H)
\tkzGetPoint{O}
\tkzDrawPoint(O)
\tkzDrawLine(C,O)
\tkzDrawLine(J,O)
\tkzLabelPoint[below right](O){$(0,0)$}
\tkzDrawLines[color=red](H,B I,A J,C)
\tkzDrawLines[color=green](A,H C,I J,B)
\tkzLabelPoint[below left](A){$(-\frac{1}{c_2x},0)$}
\tkzLabelPoint[below left](C){$(-c_1/ x,0)$}
\tkzLabelPoint[left](I){$(0,-\frac{c_1}{y})$}
\tkzLabelPoint[below left](H){$(0,-\frac{1}{y})$}
\tkzLabelPoint[above right](B){$(-1/x,0)$}
\tkzLabelPoint[right](J){$(0,-\frac{1}{c_1 c_2 y})$}
  \tkzDrawPoints(A,B,C,H,I,J)
                 \tkzLabelLine[pos=0.3, right](H,B){$r_1$}
		\tkzLabelLine[pos=0.5, right](J,C){$r_2$}
		\tkzLabelLine[pos=0.45, right](I,A){$r_3$}
		\tkzLabelLine[pos=0.4, right](I,C){$g_1$}
		\tkzLabelLine[pos=0.5, left](H,A){$g_2$}
		\tkzLabelLine[pos=0.6, right](J,B){$g_3$}
		\tkzLabelLine[pos=-0.15, above](C,O){$exc: \ \gamma=0$}
		\tkzLabelLine[pos=-0.2, right](J,O){$so: \ \beta=0$}

\end{tikzpicture}
\caption{The Pappus-Brianchon-Pascal $(9_3)_1$configuration after a projective transformation} \label{fig:Conf5}
\end{figure}

The picture, associated to the non-uniqueness factor (\ref{Q33}) is given in Figure \ref{fig:Conf5}, which is the $(9_3)_1$
 configuration after a projective transformation, which takes the $\alpha=0$ line to infinity.\footnote{One has to take into account that as the equations of the three distinguished lines are $\alpha=0, \beta=0$ and $\gamma=0$,
one of them unavoidably will be the ideal line of the projective plane, i.e. the line in the infinity (we choose $\alpha=0$). }
 
 Finally, the geometrical roles of the free parameters $c_1,c_2,x,y$ in the 
 (\ref{Q33}) expression of $Q$ are easily observed in the same Figure 
 \ref{fig:Conf5}, where the 
 associated coordinates of the points of the configuration are shown explicitly.  

\subsection{The $(16_3 12_4)$ configuration and $Q$ for the $sl, so, exc$ and $sp$
lines}

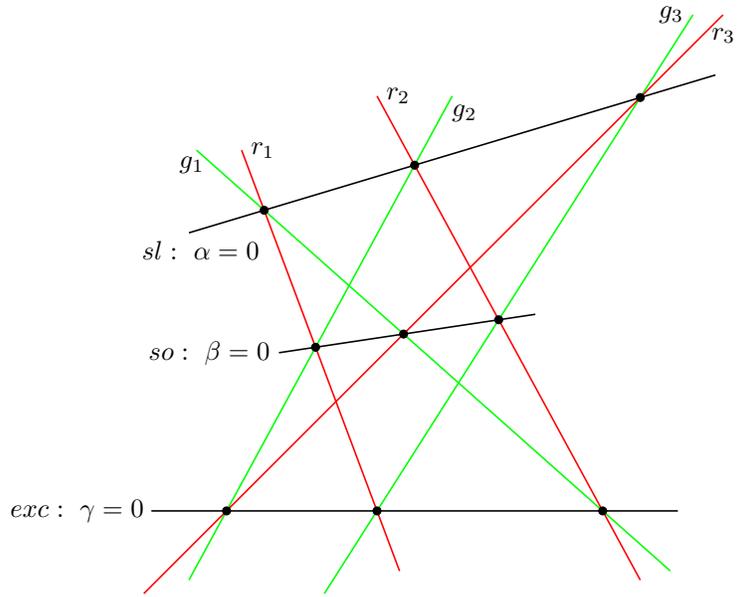
\begin{figure}
 \centering
\begin{tikzpicture} \tkzDefPoints{0/0/A,2/0/B, 5/0/C, 0.5/4/D, 2.5/4.6/E, 5.5/5.5/F}
\tkzDrawLine(A,C)
\tkzDrawLine(D,F)
\tkzDrawLines[color=green](A,E B,F)
\tkzDrawLines[color=red](C,E B,D)
\tkzDrawLines[color=red](A,F)
\tkzDrawLines[color=green](C,D)
\tkzInterLL(A,F)(C,D)
\tkzGetPoint{I}
\tkzInterLL(A,E)(B,D)
\tkzGetPoint{H}
\tkzInterLL(B,F)(C,E)
\tkzGetPoint{J}
\tkzDrawLine(H,J)
\tkzDrawPoint(H)
\tkzDrawPoint(J)
\tkzDrawPoint(I)
 \tkzDrawPoints(A,B,C,D,E,F) 
                 \tkzLabelLine[pos=-0.2, right](D,B){$r_1$}
		\tkzLabelLine[pos=-0.2, right](E,C){$r_2$}
		\tkzLabelLine[pos=-0.15, right](F,A){$r_3$}
		\tkzLabelLine[pos=-0.15, left](D,C){$g_1$}
		\tkzLabelLine[pos=-0.15, right](E,A){$g_2$}
		\tkzLabelLine[pos=-0.2, left](F,B){$g_3$}
		\tkzLabelLine[pos=-0.35, right](D,F){$sl: \ \alpha=0$}
		\tkzLabelLine[pos=-0.2, left](H,J){$so: \ \beta=0$}
		\tkzLabelLine[pos=-0.6, right](A,C){$exc: \  \gamma=0$}
\end{tikzpicture}
\caption{The Pappus-Brianchon-Pascal, or $(9_3)_1$ configuration} \label{fig:Conf2}
\end{figure}

Let's now take $4$ black lines - $sl, so, sp, exc$, and search for a
 non-uniqueness factor $Q$, which is equal to $1$ on those. 
If we take $k=4$, we happen to be dealing with the configuration  
$(16_3 12_4)$. One of its realizations, taken from \cite{GB}, is
 presented in Figure \ref{fig:Conf2}.

{\bf Proposition 4.}
{\it The configuration $(16_3 12_4)$ presented in Figure \ref{fig:Conf2} corresponds to the following non-uniqueness factor $Q$ for universal dimensions:
\begin{multline}
	Q= \frac{n \alpha + x \beta - y \gamma}{-(n+3x+3x^\prime) \alpha + x \beta - y \gamma} \cdot
	\frac{-(n+3x+3x^\prime) \alpha + x^\prime \beta - y \gamma}{n \alpha + x^\prime \beta - y \gamma} \times \\
	\frac{n \alpha + x^\prime \beta + y \gamma}{-(n+3x+3x^\prime) \alpha + x^\prime \beta + y \gamma} \cdot
	\frac{-(n+3x+3x^\prime) \alpha + x \beta + y \gamma}{n \alpha + x \beta + y \gamma}
\end{multline}
with $n, y, x$ and $x^{\prime}$ free parameters and which is equal to $1$ on each of the $sl, so, sp, exc$ distinguished lines in Vogel's plane. }
\begin{proof}
First note, that it is possible to color the lines of this particular configuration as needed: 
with the chosen coloring, demonstrated in Figure \ref{fig:Conf2}, at each of the points of the configuration three lines of different - black, red, and green - colors meet.

After labeling the lines, we track the incidence relations, which right away identify the patterns of cancellation of the factors in $Q$, or, equivalently,
 the choice of permutations in (\ref{Qx27b})-(\ref{4lineEeq}) equations; permutations, dictated by this configuration are easily identified and are as follows:

$$s(1)=2, s(2)=1, s(3)=4, s(4)=3$$
$$p(1)=4, p(2)=3, p(3)=2, p(4)=1$$
$$v(1)=3, v(2)=4, v(3)=1, v(4)=2$$

In fact, this set of permutations is the main and only information we take out of a given configuration.  
Then, having these permutations at hand, we solve the equations (\ref{4lineBeq}) - (\ref{4lineEeq}), to derive the required non-uniqueness factor.

Equations (\ref{4lineBeq})-(\ref{end1}), derived for three distinguished 
lines: $\alpha=\beta=\gamma=0$, yield

\begin{eqnarray}
	x_3= \frac{x_2}{c_2}, \,\, x_4= \frac{x_1}{c_2 k_1 k_3}  \\
	y_4=\frac{y_3}{k_3}, \,\, y_2=\frac{y_1}{k_1}          \\
	n_4=\frac{n_2 k_1}{c_2 k_1 k_3}, \,\, n_3=\frac{n_1}{k_1 c_2} \\
	k_4=\frac{1}{k_3}, \,\, k_2=\frac{1}{k_1}, \\
	c_4=\frac{1}{c_2 k_1 k_3}, \,\, c_1=c_2 k_1 k_3
\end{eqnarray}

The remaining (\ref{4-th-eq})-(\ref{end2}) equations yield 

\begin{eqnarray}
	r_3=\frac{1}{r_1}, \,\,r_2=\frac{r_1 k_3}{k_1}, \,\,  r_4=\frac{k_1}{r_1 k_3} \\
	y_1=r_1 y_3
\end{eqnarray}

and 

\begin{eqnarray}
	r_1= \pm c_2 k_1
\end{eqnarray}

The case with the plus sign, i.e. when $r_1=c_2 k_1$ leads to

\begin{eqnarray}
	n_1= k_1 n_2, \,\, x_1=k_1 x_2
\end{eqnarray}
which produces a trivial $Q$, i.e. the solution is $Q\equiv1$.

When the minus sign is taken, i.e.  $r_1=-c_2 k_1$, then:

\begin{eqnarray}
	n_2=-\frac{n_1+3x_1+3k_1 x_2}{k_1}
\end{eqnarray}

Overall, one substitutes the values of all variables, given above, expressed in terms of the remaining arbitrary variables
 $k_1, k_3, c_2, x_1, x_2, y_3, n_1$ into
 the general expression for $Q$. Absorbing some variables and defining: 
 
\begin{eqnarray}
	n=n_1, \,\, y=c_2 k_1 y_3, \,\, x=x_1, \,\, x^{\prime}=k_1 x_2
\end{eqnarray} 

we get to the final expression for the desired non-uniqueness factor $Q$ for dimensions: 

\begin{multline}
	Q= \frac{n \alpha + x \beta - y \gamma}{-(n+3x+3x^\prime) \alpha + x \beta - y \gamma} \cdot
	\frac{-(n+3x+3x^\prime) \alpha + x^\prime \beta - y \gamma}{n \alpha + x^\prime \beta - y \gamma} \times \\
	\frac{n \alpha + x^\prime \beta + y \gamma}{-(n+3x+3x^\prime) \alpha + x^\prime \beta + y \gamma} \cdot
	\frac{-(n+3x+3x^\prime) \alpha + x \beta + y \gamma}{n \alpha + x \beta + y \gamma}
\end{multline}
\end{proof}

{\bf Remark.} In case of quantum dimension the corresponding
 non-uniqueness factor $Q(x)$  appears to write almost similarly - just with the $sinh[x:$ signs inserted. 

\begin{multline}
	Q(x)= \sinh[x: \frac{n \alpha + x \beta - y \gamma}{-(n+3x+3x^\prime) \alpha + x \beta - y \gamma} \cdot
	\frac{-(n+3x+3x^\prime) \alpha + x^\prime \beta - y \gamma}{n \alpha + x^\prime \beta - y \gamma} \times \\
	\frac{n \alpha + x^\prime \beta + y \gamma}{-(n+3x+3x^\prime) \alpha + x^\prime \beta + y \gamma} \cdot
	\frac{-(n+3x+3x^\prime) \alpha + x \beta + y \gamma}{n \alpha + x \beta + y \gamma}
\end{multline}

\subsection{The $(144_3 36_{12})$ configuration and a symmetric \\ 
non-uniqueness factor $Q$}

An immediate problem, arising after the previous investigation, is the derivation of
a symmetric non-uniqueness factor $Q$, which would be equivalent to $1$ on all the
$12$ lines, obtained by the basic lines after all possible permutations of the coordinates.
The search of such a $Q$ appears to be one of a realizable $(144_3 36_{12})$ configuration in the scope of the geometrical approach.
Unfortunately, this configuration has not been studied yet, so the existence
 of a symmetric $Q$ remains an open question.

\subsection{Permutation invariance approach}

Another possible approach for searching a totally symmetric $Q$, which would be equivalent to $1$ on all lines (\ref{lines}) - (\ref{lines2}),
 is based just on its symmetry. If we find such a $Q$, which is $1$ on all three basic lines $sl, so, exc$ and is symmetric w.r.t. all permutations of the universal
  parameters, then it would be equivalent to $1$ on all permuted lines as well. 
  It is also sufficient to require its invariance w.r.t. the generating elements of the $S_3$ group of permutations, say the 
  transpositions $\alpha \leftrightarrow \beta$ and $\beta \leftrightarrow \gamma$ (in unprimed variables). 
   So in fact, a possible strategy is to take the general solution (\ref{Qx27}) for some $k$ (perhaps at least $k=12$), and require its symmetry w.r.t. these two permutations,
    which in terms of the primed parameters looks like:

\begin{eqnarray}
	\begin{vmatrix}
		\alpha^\prime \\
		\beta^\prime \\
		\gamma^\prime
	\end{vmatrix} 
	\rightarrow  \begin{vmatrix}
		\alpha^\prime \\
		3\alpha^\prime -\beta^\prime \\
		\gamma^\prime
	\end{vmatrix} \,\,\, \text{and} \,\,\,  \begin{vmatrix}
		\alpha^\prime \\
		\beta^\prime \\
		\gamma^\prime
	\end{vmatrix} 
	\rightarrow  
	\begin{vmatrix}
		\alpha^\prime +\beta^\prime+\gamma^\prime\\
		2\beta^\prime+\gamma^\prime  \\
		-3\beta^\prime-2\gamma^\prime
	\end{vmatrix},
\end{eqnarray}
repectively.

\section*{Conclusion}

We examine the problem of uniqueness of universal (quantum) dimensions of the simple Lie algebras.
In other words, we search for non-uniqueness factors, i.e. non-trivial expressions, yielding $1$ at the points
from Vogel's table. They obviously preserve the values of a universal formula at the corresponding point, when
multiplying them by the latter.

We derive several non-uniqueness factors, particularly one, which is equivalent to 
$1$, when restricting it to each of the
 three - $sl, so, exc$ - distinguished lines in the Vogel's plane, as well as another one, which
is equivalent to $1$ on each of the four distinguished lines - $sl, so, sp, exc$.
 
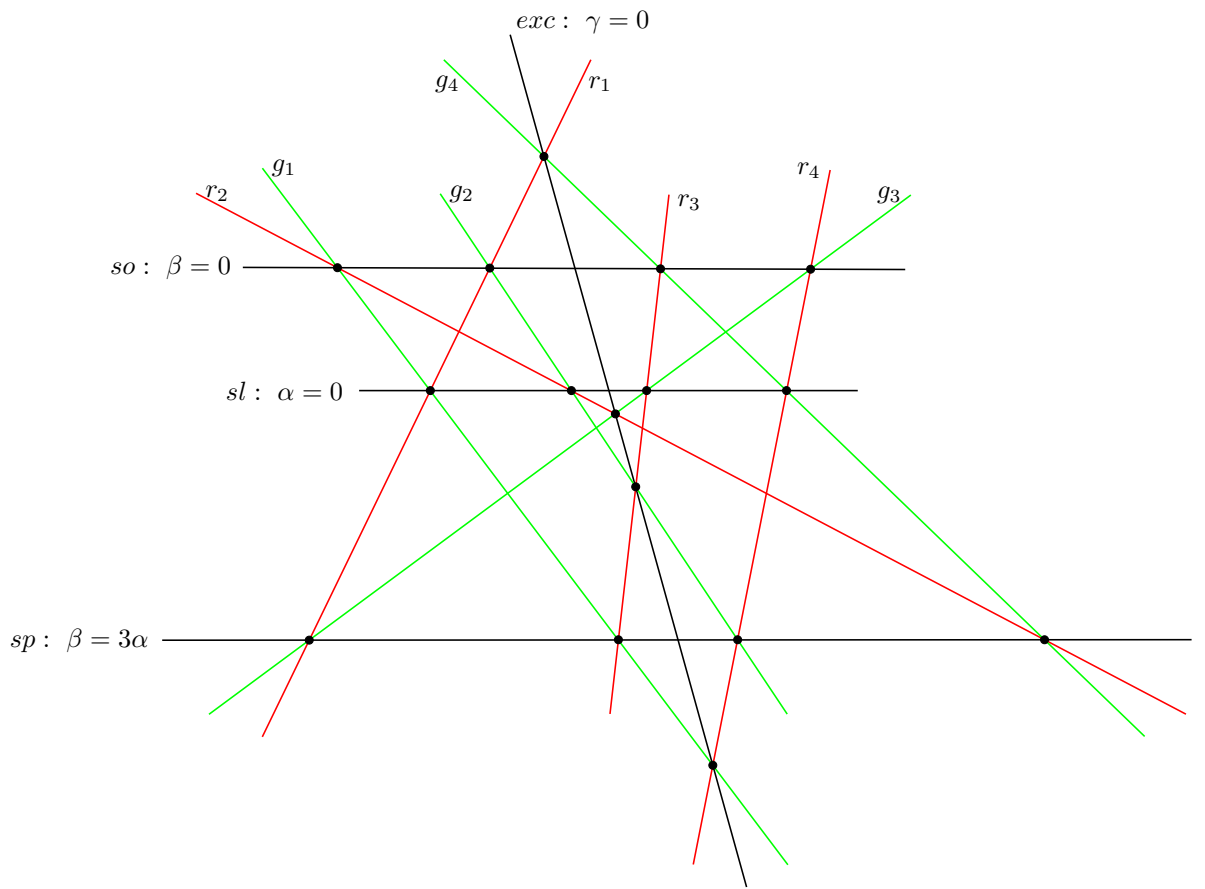
\begin{figure}
\centering
	\begin{tikzpicture} \tkzDefPoints{5/5/A,6/5/B, 4.635/8.115/C, 
			5.585/4.69/D, 5.855/3.72/E, 1.515/1.68/F, 11.29/1.685/G, 3.915/6.63/H,
			6.185/6.62/I, 3.125/5/L, 7.86/5/M, 1.89/6.635/N, 8.18/6.615/O,
			5.625/1.685/J, 7.21/1.685/K, 6.88/0.015/P}
		\tkzDrawLine(C,P)
		\tkzDrawLine(N,O)
		\tkzDrawLine(L,M)
		\tkzDrawLine(F,G)
		\tkzDrawLines[color=red](N,G I,J F,C O,P)
		\tkzDrawLines[color=green](N,P C,G F,O H,K)
		\tkzLabelLine[pos=-0.2, right](N,P){$g_1$}
		\tkzLabelLine[pos=-0.2, right](H,K){$g_2$}
		\tkzLabelLine[pos=-0.15, left](C,G){$g_4$}
		\tkzLabelLine[pos=1.2, left](F,O){$g_3$}
		\tkzLabelLine[pos=-0.15, right](C,F){$r_1$}
		\tkzLabelLine[pos=-0.2, right](N,G){$r_2$}
		\tkzLabelLine[pos=-0.18, right](I,J){$r_3$}
		\tkzLabelLine[pos=-0.2, left](O,P){$r_4$}
		\tkzLabelLine[pos=-0.22, right](C,P){$exc: \  \gamma=0$}
		\tkzLabelLine[pos=-0.5, right](N,O){$so: \  \beta=0$}
		\tkzLabelLine[pos=-0.6, right](L,M){$sl: \  \alpha=0$}
		\tkzLabelLine[pos=-0.42, right](F,G){$sp: \  \beta=3\alpha$}
				\tkzDrawPoints(A,B,C,D,E,F,G,H,I,L,M,N,O,J,K,P)
			\end{tikzpicture}
	\caption{The $16_3 12_4$ configuration} \label{fig:Conf6}
\end{figure}

 Derivation of the latter has been carried out,
 using the remarkable connection of the problem of uniqueness with the theory of 
 configurations of points and lines, which we establish by setting up a geometric interpretation
  of the non-uniqueness factors; 
  it particularly corresponds to the $(16_3 12_4)$ configuration, presented in Figure (\ref{fig:Conf6}).
  
  More interestingly, this geometric interpretation of the uniqueness problem reveals the one-to-one
  correspondence of the non-uniqueness factor, which yields $1$ on three distinguished lines, with
  the famous $(9_3)_1$ Pappus-Brianchon-Pascal configuration, known since the $4$th century AD!
  
 However, the non-uniqueness factors, presented above, will violate 
 the reasonable outputs of universal formulae, when considering
 it after some permutation of the coordinates \cite{AM, X2kng}. 
 To preserve the outputs at the points with permuted coordinates, one needs to have a completely symmetric factor. 
 In the scope of the theory of configurations, the required non-uniqueness factor 
 will probably be connected with the $(144_3 36_{12})$ configuration, which is yet to be studied.

The problem of uniqueness deserves a further investigation. 
Its connection with the classical problem of configurations of points and lines 
is intriguing and seems to have a potential of possible influence in both directions.

\section*{Acknowledgments}

We are indebted to the referee of our paper \cite{X2kng} for a
 question which is partially answered by the present investigation. MA is
  grateful to the organizers of RDP Online Workshop on Mathematical Physics
(December 5-6, 2020) for invitation.

The work of MA was fulfilled within the Regional Doctoral Program on 
Theoretical and Experimental Particle Physics
sponsored by VolkswagenStiftung. 
The work of MA and RM is partially supported by the Science Committee of
 the Ministry of Science 
and Education of the Republic of Armenia under contract  20AA-1C008.

\section*{Data availability statement}
Data sharing not applicable to this article as no datasets were generated or analyzed during the current study.

\begin{table}[ht] 
	\caption{Vogel's parameters for simple Lie algebras and the distinguished lines}
	\centering
	\begin{tabular}{|r|r|r|r|r|r|} 
		\hline Algebra/Parameters & $\alpha$ &$\beta$  &$\gamma$  & $t$ & Line \\ 
		\hline  $sl(N)$  & -2 & 2 & $N$ & $N$ & $\alpha+\beta=0$ \\ 
		\hline $so(N)$ & -2  & 4 & $N-4$ & $N-2$ & $ 2\alpha+\beta=0$ \\ 
		\hline $sp(2N)$ & -2  & 1 & $N+2$ & $N+1$ & $ \alpha+2\beta=0$ \\ 
		\hline $exc(n)$ & $-2$ & $2n+4$  & $n+4$ & $3n+6$ & $\gamma=2(\alpha+\beta)$\\ 
		\hline 
	\end{tabular}
	
	{On the $exc$ line $n=-2/3,0,1,2,4,8$ for $G_2, so(8), F_4, E_6, E_7,E_8 $, 
		respectively.} \label{tab:V2}
\end{table}


\newpage

\begin{thebibliography}{99}

\bibitem{V0}
P. Vogel,  The Universal Lie algebra. Preprint (1999), https://webusers.imj-prg.fr/\~{}pierre.vogel/grenoble-99b.pdf

\bibitem{V}
P.Vogel,   Algebraic structures on modules of diagrams. Preprint (1995),  www.math.jussieu.fr/\~{}vogel/diagrams.pdf, J. Pure Appl. Algebra {\bf 215} (2011), no. 6, 1292-1339. 


\bibitem{Del}
P. Deligne, La s\'erie exceptionnelle des groupes de Lie, C. R. Acad. Sci. Paris, S\'erie I {\bf 322} (1996), 321-326.

\bibitem {DM}
P. Deligne and R. de Man, La s\'erie exceptionnelle des groupes de Lie II, C. R. Acad. Sci. Paris, S\'erie I  {\bf 323} (1996), 577-582. 


\bibitem{Cohen}
A. M.Cohen and  R. de Man,   Computational evidence for Deligne's conjecture regarding exceptional Lie
groups, Comptes Rendus de l'Académie des Sciences, Série 1, Mathématique, (1996) 322(5), 427-432	

\bibitem{W3}
B.Westbury, Invariant tensors and diagrams, Proceedings of the Tenth Oporto Meeting on Geometry, Topology and Physics (2001). Vol. 18. October, suppl. 2003, pp. 49-82.

\bibitem{LM1} 
J.M. Landsberg and  L.Manivel,  A universal dimension formula for complex simple Lie algebras. Adv. Math. {\bf 201} (2006), 379-407

	\bibitem{MV}
R.L. Mkrtchyan  and A.P.Veselov,  Universality in Chern-Simons theory. JHEP08 (2012) 153, arxiv:1203.0766.	

\bibitem{M16QD}
R.L.Mkrtchyan, On Universal Quantum Dimensions, arxiv:1610.09910, Nuclear Physics B921,  2017, pp. 236-249, 


 \bibitem{AM}
 M.Y. Avetisyan and R.L. Mkrtchyan,  $X_2$ Series of Universal Quantum Dimensions, arXiv:1812.07914, J. Phys. A: Math. Theor. Volume 53, Number 4, 045202,  https://doi.org/10.1088/1751-8121/ab5f4d 

\bibitem{X2kng}
M.Y. Avetisyan and R.L. Mkrtchyan, On $(ad)^n(X_2)^k$ series of universal quantum dimensions, arXiv:1909.02076, J. Math. Phys. 61, 101701 (2020), https://doi.org/10.1063/5.0007028


\bibitem{AMLR}
M.Y. Avetisyan and R.L. Mkrtchyan, On linear resolvability of universal quantum dimensions, arXiv:2101.08780.


\bibitem{GKF}
Hilbert, David; Cohn-Vossen, Stephan (1952), Geometry and the Imagination (2nd ed.), Chelsea,  ISBN 0-8284-1087-9.

\bibitem{GB}
Grünbaum, Branko (2009), Configurations of Points and Lines, Graduate Studies in Mathematics, 103, American Mathematical Society, ISBN 978-0-8218-4308-6.




\end{thebibliography}
\end{document}